\documentclass{amsart}

\newtheorem{theorem}{Theorem}
\newcommand{\bt}{\begin{theorem}}
\newcommand{\et}{\end{theorem}}
\newtheorem{lemma}{Lemma}
\newcommand{\bl}{\begin{lemma}}
\newcommand{\el}{\end{lemma}}
\newtheorem{corollary}{Corollary}
\newcommand{\bc}{\begin{corollary}}
\newcommand{\ec}{\end{corollary}}
\newtheorem{problem}{Problem}
\newcommand{\bprob}{\begin{problem}}
\newcommand{\eprob}{\end{problem}}
\newcommand{\beq}{\begin{equation}}
\newcommand{\eeq}{\end{equation}}
\newcommand{\benum}{\begin{enumerate}}
\newcommand{\eenum}{\end{enumerate}}
\newcommand{\N}{\ensuremath{ \mathbf N }}

\newcommand{\R}{\ensuremath{\mathbf R}}
\newcommand{\C}{\ensuremath{\mathbf C}}
\newcommand{\mca}{\ensuremath{ \mathcal A}}
\newcommand{\mcb}{\ensuremath{ \mathcal B}}

\newcommand{\mce}{\ensuremath{ \mathcal E}}
\newcommand{\mcf}{\ensuremath{ \mathcal F}}

\newcommand{\mcj}{\ensuremath{ \mathcal J}}

\newcommand{\mcx}{\ensuremath{ \mathcal X}}

\DeclareMathOperator{\card}{\text{card}}

\newcommand{\bmat}{\left(\begin{matrix}}
\newcommand{\emat}{\end{matrix}\right)}

\DeclareMathOperator{\qqand}{\qquad\text{and}\qquad}

\title{Weighted real Egyptian numbers}
\author{Melvyn B. Nathanson}
\address{Department of Mathematics\\
Lehman College (CUNY)\\
Bronx, NY 10468}
\email{melvyn.nathanson@lehman.cuny.edu}
\date{\today}

\thanks{Supported in part by a grant from the PSC-CUNY Research Award Program.}

\begin{document}

\maketitle

\begin{abstract}
Let $\mca = (A_1,\ldots, A_n)$ be a  sequence of nonempty finite sets 
of positive real numbers, 
and let $\mcb = (B_1,\ldots, B_n)$ be a sequence of infinite discrete 
sets of positive real numbers.  
A  \emph{weighted real Egyptian number with numerators \mca\ 
and denominators \mcb} 
is a real number $c$ that can be represented in the form
\[
c = \sum_{i=1}^n \frac{a_i}{b_i}
\]
with $a_i \in A_i$ and $b_i \in B_i$ for $i \in \{1,\ldots, n\}$.   
In this paper, classical results of Sierpi\' nski for Egyptian fractions 
are extended to the set of weighted real Egyptian numbers.  
\end{abstract}

\section{Weighted Egyptian numbers}

Let $\N = \{1,2,3,\ldots\}$ denote the set of positive integers.

An \emph{Egyptian fraction of length $n$} is a rational number that can be represented
as the sum of $n$ pairwise distinct unit fractions, that is, a rational number of the form
\[
\sum_{i=1}^n \frac{1}{b_i}
\]
for some $n$-tuple $(b_1,\ldots, b_n)$ of pairwise distinct positive integers. 
Deleting the requirement that the denominators be pairwise distinct, 
we define an \emph{Egyptian number of length $n$} as a rational number 
that is the sum of $n$ unit fractions, that is, a rational number of the form
\[
\sum_{i=1}^n \frac{1}{b_i}
\]
for some $n$-tuple $(b_1,\ldots, b_n)$ of positive integers.  
Because $1/b = 1/2b + 1/2b$, an Egyptian number of length at most $n$ 
is also an Egyptian number of length $n$.  

Repeated use of the elementary identities
\[
\frac{2}{2k} =\frac{1}{k} = \frac{1}{k+1} + \frac{1}{k(k+1)}
\]
and
\[
\frac{2}{2k+1} = \frac{1}{k+1} + \frac{1}{(k+1)(2k+1)}
\]
allows us to write  an Egyptian number of length $n$ 
as an Egyptian fraction of length $n$,   
and also to write an Egyptian fraction of length $n$ 
as an Egyptian fraction  of length $n'$ for every $n' \geq n$.

Richard K. Guy's book 
\emph{Unsolved Problems in Number Theory}~\cite[pp. 252--262]{guy04}
contains an ample bibliography and many open questions about Egyptian fractions.

There is a natural extension of Egyptian numbers from the the set 
of positive rational numbers to the set of positive real numbers.  
Let $A$ be a finite set of positive real numbers, 
and let $B$ be an infinite discrete set of positive real numbers.  
(The set $B$  is \emph{discrete} if $B \cap X$ is finite 
for every bounded set $X$.)
We consider ``unit fractions'' of the form $1/b$ with $b \in B$, 
and finite sums of these unit fractions with weights $a \in A$.  
This gives  real numbers of the form $\sum_{i=1}^n a_i/b_i$.

More generally, let  $\mca = (A_1,\ldots, A_n)$ be a  sequence of nonempty finite sets 
of positive real numbers, 
and let $\mcb = (B_1,\ldots, B_n)$ be a sequence of infinite discrete 
sets of positive real numbers.  
A  \emph{weighted real Egyptian number with numerators \mca\ 
and denominators \mcb} 
is a real number $c$ that can be represented in the form
\[
c = \sum_{i=1}^n \frac{a_i}{b_i}
\]
for some 
\[
(a_1,\ldots, a_n,b_1,\ldots, b_n) \in 
A_1 \times \cdots \times A_n \times B_1 \times \cdots \times B_n.
\]  
Let 
\[
\mce(\mca,\mcb) = \left\{ \sum_{i=1}^n \frac{a_i}{b_i} : a_i \in A_i \text{ and } b_i \in B_i 
\text{ for } i \in \{1,\ldots, n\}  \right\}
\]
 be the set of all weighted real Egyptian numbers with numerators \mca\ 
and denominators \mcb.   The set $\mce(\mca,\mcb)$ is a set of positive real numbers. 

For all $c \in \R$, we  define the \emph{representation function}
\[
r_{\mca,\mcb}(c) = \card\left( (a_1,\ldots, a_n,b_1,\ldots, b_n) \in 
A_1 \times \cdots \times A_n \times B_1 \times \cdots \times B_n :
 \sum_{i=1}^n \frac{a_i}{b_i} = c \right). 
\] 

The purpose of this note is to show that the topological results 
about Egyptian numbers in 
Sierpi{\' n}ski's classic paper~\cite{sier56}, 
``Sur  les decompositions de nombres rationnels en fractions primaires''
extend to weighted Egyptian numbers.

Note that an Egyptian number of length $n$ 
is a weighted real Egyptian number 
with numerators $\mca = (\{ 1\}, \ldots, \{ 1\})$ 
and denominators $\mcb = (\N,\ldots, \N)$.  
Conversely, for all $a,b \in \N$, we have 
\[
\frac{a}{b} = \underbrace{\frac{1}{b} + \cdots +  \frac{1}{b}}_{\text{$a$ summands}}.
\]
Thus, every weighted  real Egyptian number 
with numerators $\mca = (A_1,\ldots, A_n)$ such that $A_i$ is a finite set 
of positive integers for $i \in \{1,\ldots, n\}$,
and with denominators $\mcb = (\N,\ldots, \N)$,    
is an Egyptian number of length at most $\sum_{i=1}^n \max(A_i)$.

\bt           \label{sierpinski:theorem:fg-new}
Let $A_1,\ldots, A_n$ be nonempty finite sets of positive real numbers, and let 
$B_1,\ldots, B_n$ be infinite discrete sets of  positive real  numbers.
Let       
\beq                       \label{sierpinski:2n-tuple-new}
\left( (a_{m,1},\ldots, a_{m,n}, b_{m,1},\ldots, b_{m,n}) \right)_{m \in \N} 
\eeq
be an infinite sequence of pairwise distinct $2n$-tuples  
in $A_1\times \cdots \times A_n \times B_1\times \cdots \times B_n$,  
that is, 
\[
 (a_{m,1},\ldots, a_{m,n}, b_{m,1},\ldots, b_{m,n}) 
 =  (a_{m',1},\ldots, a_{m',n}, b_{m',1},\ldots, b_{m',n})
\]
if and only if  $m = m'$.
For $m \in \N$, let 
\[
c_m = \sum_{i=1}^n \frac{ a_{m,i} } {b_{m,i} } \in \mce(\mca,\mcb).  
\]   
The sequence $(c_m)_{m\in \N}$ contains a strictly decreasing subsequence.
\et

Equivalently, there exists a strictly increasing sequence $\left( m_j \right)_{j =1}^{\infty}$ 
of positive integers such that 
\[
c_{m_j } >  c_{m_{j+1} } > 0
\]
for all $j \in \N$.

\begin{proof}
For $i  \in \{ 1,2,\ldots, n\}$, let 
\[
B_{0,i} = \{b_{m,i}: m = 1,2,3,\ldots\}
\]
where $b_{m,i}$ is the $(n+i)$-th coordinate of the $m$th 
$2n$-tuple in the sequence~\eqref{sierpinski:2n-tuple-new}.
We have $B_{0,i}   \subseteq B_i$ and 
\[
 (a_{m,1},\ldots, a_{m,n}, b_{m,1},\ldots, b_{m,n}) 
 \in A_1\times \cdots \times A_n \times B_{0,1} \times \cdots \times B_{0,n}
\]
for all $m \in \N$.
If the set $B_{0,i} $ is finite for all $i = 1,\ldots, n$, then the set 
$ A_1\times \cdots \times A_n \times B_{0,1} \times \cdots \times B_{0,n} $ 
is finite.  This implies that the sequence~\eqref{sierpinski:2n-tuple-new} 
is finite, which is absurd.  
Therefore, $B_{0,i} $ is infinite for some $i \in \{1,\ldots, n\}$.  
Without loss of generality, we can assume that $i = 1$ 
and $B_{0,1} $ is infinite.  

Because $B_{0,1} $ is contained in the discrete set $B_1$, 
there is a strictly increasing sequence 
of positive integers $(m_{j,1})_{j=1}^{\infty}$ such that 
\[
\lim_{j \rightarrow \infty} b_{m_{j,1},1} = \infty.
\]
Let $k \in \{1,\ldots, n\}$, and let $(m_{j,k})_{j=1}^{\infty}$ be 
a strictly increasing sequence of positive integers  such that 
\[
\lim_{j \rightarrow \infty} b_{m_{j,k},i} = \infty
\]
for $i \in \{1,\ldots, k \}$.  
If $k \leq n-1$, then, for $i \in \{ k+1,k+2,\ldots, n\}$, we consider the set 
\[
B_{k,i} = \{ b_{m_{j,k},i} : j \in \N \}.
\]
Suppose that the set $B_{k,i}$ is infinite for some $i \in \{ k+1,k+2,\ldots, n\}$.   
Without loss of generality, we can assume that $i = k+1$.  Because 
$B_{k,k+1}$ is an infinite subset of the discrete set $B_{k+1}$, 
the sequence $(m_{j,k})_{i=1}^{\infty}$ contains a strictly increasing subsequence
$(m_{j,k+1})_{j=1}^{\infty}$ such that 
\[
\lim_{j\rightarrow \infty} b_{m_{j,k+1},k+1} = \infty.
\]
It follows that 
\[
\lim_{j\rightarrow \infty} b_{m_{j,k+1}, i} = \infty      
\]
for all $i \in \{1, 2,\ldots, k, k+1\}$.  
Continuing inductively, we obtain an integer $s \in \{1,2,\ldots, n\}$ 
and a strictly increasing sequence of positive integers $(m_{j,s})_{j =1}^{\infty}$ such that 
\beq         \label{sierpinski:g-limit-new}
\lim_{j\rightarrow \infty} b_{m_{j,s}, i} = \infty 
\eeq
for all $i \in \{1, 2,\ldots, s\}$, 
and the sets
\[
B_{s,i} =  \{ b_{m_{j,s},i} :  j \in \N \}
\]
are finite for all $i \in \{s+1,\ldots, n\}$.  

The sets $A_1,\ldots, A_n$ and $B_{s,s+1},\ldots, B_{s,n}$ are finite. 
Therefore, the set of $(2n-s)$-tuples
\[
A_1 \times \cdots \times A_n \times B_{s,s+1} \times \cdots \times B_{s,n}
\]
is finite.  By the pigeonhole principle, there exists a $(2n-s)$-tuple 
\[
\left(a^*_1,\ldots, a^*_n, b^*_{s+1}, \ldots, b^*_{n}\right) 
\in A_1 \times \cdots \times A_n \times B_{s,s+1} \times \cdots \times B_{s,n}
\]
and a strictly increasing subsequence $(m_{j,s+1})_{t=1}^{\infty}$ of the sequence 
$(m_{j,s})_{j=1}^{\infty}$ such that 
\[
\left( a_{m_{j,s+1},1} ,\ldots, a_{m_{j,s+1},n} , b_{m_{j,s+1},s+1} , \ldots, b_{m_{j,s+1},n}\right)
= \left(a^*_1,\ldots, a^*_n, b^*_{s+1}, \ldots, b^*_{n}\right)
\]
for all $j \in \N$.     
It follows that, for all $j \in \N$, 
\begin{align*}
c_{ m_{j,s+1}}
& = \sum_{i=1}^s \frac{a^*_i }{ b_{m_{j,s+1},i}} 
+ \sum_{i=s+1}^n \frac{ a^*_i}{ b^*_i } 
 = \sum_{i=1}^s \frac{a^*_i }{ b_{m_{j,s+1},i} } + c^*_0
\end{align*}
where
\[
c^*_0 = \sum_{i=s+1}^n \frac{ a^*_i}{ b^*_i } \geq 0.
\]
Note that $c^* > 0$ if $s < n$ and $c^* = 0$ if $s=n$.

The limit condition~\eqref{sierpinski:g-limit-new} implies that 
there exists a strictly increasing sequence of positive integers 
$( m_{j,s+2} )_{j=1}^{\infty}$ such that 
\[
 b_{m_{j,s+2}, i}   <  b_{ m_{j+1,s+2}, i} 
\]
for all $i \in \{1,\ldots, s\}$ and for all $j \in \N$.
Let 
\[
m_j = m_{j,s+2}
\]
for $j \in \N$.       
We have
\[
b_{m_j, i}  < b_{ m_{j+1}, i} 
\]
for all $i \in \{1,\ldots, s\}$, and so
\[
c_{ m_j}
 = \sum_{i=1}^s \frac{a^*_i}{b_{m_j,i}  } + c^*_0 
 >  \sum_{i=1}^s \frac{f_i(a^*_i )}{b_{m_{j+1}, i}} + c^*_0 
 = c_{ m_{j+1}} > 0
\]
for all $j \in \N$.  
This completes the proof.
\end{proof}

\bc      
If $\mca = (A_1,\ldots, A_n)$ is a  sequence of nonempty finite sets 
of positive real numbers 
and $\mcb = (B_1,\ldots, B_n)$ is a sequence of infinite discrete 
sets of positive real numbers, then   
\[
r_{\mca,\mcb}(c) < \infty
\]
for all $c \in \R$.
\ec

\begin{proof}
Because $\mce(\mca,\mcb)$ is a set of positive real numbers, 
we have $r_{\mca,\mcb}(c) = 0$ for all $c \leq 0$.

If $r_{\mca,\mcb}(c) = \infty$ for some $c> 0$, 
then there exists an infinite sequence of pairwise distinct $2n$-tuples  
of the form~\eqref{sierpinski:2n-tuple-new} 
such that $c_m = c$ for all $m \in \N$, 
and the constant sequence $(c_m)_{m\in \N}$ contains no strictly decreasing subsequence.  
This is impossible by Theorem~\ref{sierpinski:theorem:fg-new}.
\end{proof}

\bc        \label{sierpinski:corollary:delta}
For every $c \in  \R$ there exists $\delta = \delta(c) > 0$ 
such that $(c-\delta, c) \cap \mce(\mca,\mcb) = \emptyset$.
\ec

\begin{proof}
Let $c \in  \R$. 
If, for every positive integer $m$, there exists 
\[
c_m \in \left(c - \frac{1}{m},c\right) \cap \mce(\mca,\mcb),
\]
then the sequence $(c_m)_{m\in \N}$  contains a strictly increasing subsequence, 
and this subsequence contains no strictly decreasing subsequence.  
This is impossible by Theorem~\ref{sierpinski:theorem:fg-new}.
 Therefore, there exists  $m \in \N$ such that $\delta = 1/m > 0$
satisfies the condition $\left(c - \delta, c \right) \cap \mce(\mca,\mcb) = \emptyset$.  
\end{proof}

\bc
The set $\mce(\mca,\mcb)$ is nowhere dense.  
\ec

\begin{proof}
Let  $\overline{\mce(\mca,\mcb)}$ denote the closure of $\mce(\mca,\mcb)$,
and let $U$ be a nonempty open set in \R.
If $U \cap \overline{\mce(\mca,\mcb)} \neq \emptyset$, 
then there exists $c \in U \cap \mce(\mca,\mcb)$.  
By Corollary~\ref{sierpinski:corollary:delta}, there exists $\delta > 0$ such 
that $(c - \delta, c) \cap \mce(\mca,\mcb) = \emptyset$, and so 
 $U \not\subseteq \overline{\mce(\mca,\mcb)}$.  
 It follows that the set  $\mce(\mca,\mcb)$ of weighted real Egyptian numbers 
 is nowhere dense.
\end{proof}

\section{Signed  weighted Egyptian numbers}

\emph{Notation.}. 
Let $j_1, \ldots, j_s \in \N$.  We write 
\[
 (j_1, \ldots, j_s) \preceq (1,\ldots, n\}
\]
if $1 \leq j_1 < j_2 < \cdots < j_s \leq n$.
For $s \in \{1,\ldots, n-1\}$ and 
\[
J = (j_1, \ldots, j_s) \preceq (1,\ldots, n\}
\]
let  
\[
L = (1,\ldots, n) \setminus J = (\ell_1, \ldots, \ell_{n-s}) 
\]
be the strictly increasing $(n - s)$-tuple obtained by deleting the integers 
$j_1,\ldots, j_s$ from $(1,\ldots, n)$.     
To the $n$-tuple of sets $\mca = \left( A_1, A_2, \ldots, A_n \right)$, 
we associate the $s$-tuple of sets 
\[
\mca_J = \left(A_{j_1}, A_{j_2},\ldots, A_{j_s} \right).     
\]
and the $(n-s)$-tuple of sets 
\[
\mca_L = \left(A_{\ell_1}, A_{\ell_2},\ldots, A_{\ell_{n-s}} \right).     
\]
For example, the 2-tuple  
\[
J = (3,5) \preceq (1,2,3,4,5,6)
\]
 and the 4-tuple 
\[
L = (1,2,3,4,5,6) \setminus (3,5) = (1,2,4,6)
\]
determine the set sequences   $\mca_J = (A_3, A_5)$ 
and $\mca_L = (A_1, A_2, A_4, A_6)$.

Let $\mca = (A_1,\ldots, A_n)$ be a  sequence of nonempty finite sets 
of positive real numbers, 
and let $\mcb = (B_1,\ldots, B_n)$ be a sequence of infinite discrete 
sets of positive real numbers.  
A  \emph{signed weighted real Egyptian number with numerators \mca\ 
and denominators \mcb} 
is a real number $c$ that can be represented in the form
\beq              \label{sierpinski:3n-tuple=c}
c = \sum_{i=1}^n \frac{\varepsilon_i a_i}{b_i}
\eeq
for some $3n$-tuple 
\begin{align}               \label{sierpinski:3n-tuple}
(a_1,  \ldots,  & a_n,b_1,\ldots, b_n, \varepsilon_1,\ldots, \varepsilon_n) \\
&  \in 
A_1 \times \cdots \times A_n \times B_1 \times \cdots \times B_n \times \{ 1, -1\}^n.
\nonumber
\end{align}
Let 
\[
\mce^{\pm}(\mca,\mcb) = \left\{ \sum_{i=1}^n \frac{\varepsilon_ia_i}{b_i} : 
a_i \in A_i, \ b_i \in B_i, \text{ and } \varepsilon_i \in \{ 1, -1 \} 
\text{ for } i \in \{1,\ldots, n\}  \right\}
\]
 be the set of all signed weighted Egyptian numbers with numerators \mca\ 
and denominators \mcb.   
For all $c \in \R$,  
the \emph{representation function} $ r^{\pm}_{\mca,\mcb}(c) $ 
counts the number of $3n$-tuples of the form~\eqref{sierpinski:3n-tuple}
that satisfy equation~\eqref{sierpinski:3n-tuple=c}.
We have $ r^{\pm}_{\mca,\mcb}(c) \geq 1$  if and only if $c \in \mce^{\pm}(\mca,\mcb) $.

The proofs in this section are simple modifications of proofs in~\cite{sier56}.

\bt              \label{sierpinski:theorem:reps1-2}
Let $\mca = (A_1,\ldots, A_n)$ be a  sequence of nonempty finite sets 
of positive real numbers, 
and let $\mcb = (B_1,\ldots, B_n)$ be a sequence of infinite discrete 
sets of positive real numbers.  
If $n=1$, then 
\[
r^{\pm}_{\mca,\mcb}(c) < \infty
\]
for all $c \in \R$.
If $n=2$, then 
\[
r^{\pm}_{\mca,\mcb}(c) < \infty
\]
for all $c \in \R\setminus \{0\}$, but it is possible that $r^{\pm}_{\mca,\mcb}(0) = \infty$.

Let $n \geq 3$.  Let $s \in \{2,3,\ldots, n-1)$,  
$J = (j_1,\ldots, j_s) \preceq (1,\ldots, n)$, and  $L = (1,\ldots, n) \setminus J$.
If $r^{\pm}_{\mca_J,\mcb_J}(0) = \infty$, then $r^{\pm}_{\mca,\mcb}(c) = \infty$ 
for all $c \in \mce^{\pm}(\mca_L,\mcb_L)$.     
\et

\begin{proof} 
If $n=1$,  $\mca = (A_1)$, and $\mcb = (B_1)$, then 
\[
\mce^{\pm}(\mca,\mcb) = \left\{  \frac{\varepsilon_1 a_1}{b_1} : 
a_1 \in A_1, \ b_1 \in B_1, \text{ and } \varepsilon_1 \in \{ 1, -1 \}  \right\}
\]
is a set of nonzero numbers, and so  $r^{\pm}_{\mca,\mcb}(0) = 0$ .  

Let $c \in \R\setminus \{0\}$.  If $r^{\pm}_{\mca,\mcb}(c) \geq 1$, then 
$c = \varepsilon_1 a_1/b_1$ for some $a_1 \in A_1, \ b_1 \in B_1, \varepsilon_1 \in \{ 1, -1 \}$.
If $c > 0$, then $\varepsilon_1 = 1$.  If $c < 0$, then $\varepsilon_1 = -1$.  
For each $a_1 \in A_1$ there is at most one $b_1 \in B_1$ such that 
$c = \varepsilon_1 a_1/b_1$, and so $r^{\pm}_{\mca,\mcb}(c) \leq |A_1| < \infty$.

Let $n=2$.  Suppose that $\mca = (A_1, A_2)$ and $\mcb = (B_1, B_2)$. 
Let $A = A_1 \cap A_2$ and $B = B_1 \cap B_2$.  
If $A$  is nonempty and $B$ is infinite, then for all $a\in A$ and $b \in B$  
we have
\[
(a,a,b,b,1,-1) \in A_1 \times A_2 \times B_1 \times B_2 \times \{1,-1\}^2 
\] 
and 
\[
0 = \frac{a}{b} + \frac{(-a)}{b}
\]
 and so  $r^{\pm}_{\mca,\mcb}(0) = \infty$.  

Let $c \in \R\setminus \{0\}$. 
Let $a^* = \max(A_1 \cup A_2)$.  
If 
\beq          \label{sierpinski:rep-c-2}
c =  \frac{\varepsilon_1 a_1}{b_1} +  \frac{\varepsilon_2 a_2}{b_2}
\eeq
is a representation of $c$ in $\mce^{\pm}(\mca,\mcb)$, then  
\[
|c| \leq \frac{ a_1}{b_1} +  \frac{ a_2}{b_2} 
\leq a^* \left( \frac{ 1}{b_1} +  \frac{ 1}{b_2}  \right) 
\leq \frac{2a^*}{ \min(b_1, b_2)} 
\]
and so 
\[
0 <  \min(b_1, b_2) \leq \frac{na^*}{|c|} 
\]
Because  the sets $B_1$ and $B_2$ are discrete, the sets  
\[
\tilde{B_i} =  \left\{b_i \in B_i:b_i \leq \frac{na^*}{|c|}  \right\}
\]
are finite for $i=1$ and 2, and so the set of fractions
\[
\mcf = \bigcup_{i=1}^2 \left\{ \frac{\varepsilon_i a_i}{b_i} : a_i \in A_i, b_i \in \tilde{B_i}, \varepsilon_i \in \{1,-1\}    \right\}
\] 
is also finite.  
Every representation of $c$ of the form~\eqref{sierpinski:rep-c-2} 
must include at least one fraction in the set \mcf, and this fraction uniquely 
determines the other fraction in the representation~\eqref{sierpinski:rep-c-2}.
Therefore, $r^{\pm}_{\mca,\mcb}(c) < \infty$ for $c \neq 0$.

The statement for $n \geq 3$ follows immediately from the observation 
that if $J \preceq (1,\ldots, n)$ and  $L = (1,\ldots, n) \setminus J$, 
then 
\[
\mce^{\pm}(\mca_J,\mcb_J) + \mce^{\pm}(\mca_L,\mcb_L) = \mce^{\pm}(\mca,\mcb).
\]
This completes the proof. 
\end{proof}

\bt                       \label{sierpinski:theorem:reps}
Let $\mca = (A_1,\ldots, A_n)$ be a  sequence of nonempty finite sets 
of positive real numbers, 
and let $\mcb = (B_1,\ldots, B_n)$ be a sequence of infinite discrete 
sets of positive real numbers.  
Let
\[
\mcj(\mca,\mcb)
=  \bigcup_{s=1}^{n-2}  \quad   
\bigcup_{ \substack{J_s =  (1_1, \ldots, j_s) \\ \preceq (1,\ldots, n\}}}
\mce^{\pm}(\mca_{J_s}, \mcb_{J_s}).
\]
For all $c \in \R \setminus \mcj(\mca,\mcb)$, 
\[
r^{\pm}_{\mca,\mcb}(c) < \infty.
\]
\et

\begin{proof}
The sets $A_1,\ldots, A_n$ are nonempty and finite.  Let 
\[
a^* = \max\left( \bigcup_{i=1}^n A_i \right).
\]
Let $c \in \R\setminus \{0\}$.  If
\[
c = \sum_{i=1}^n \frac{\varepsilon_i a_i}{b_i} 
\]
is a representation of $c$ in $\mce^{\pm}(\mca,\mcb)$, then  
\[
|c| \leq  \sum_{i=1}^n \frac{a_i}{b_i} \leq a^*\sum_{i=1}^n \frac{1}{b_i} 
\leq \frac{na^*}{ \min\{b_1,\ldots, b_n\}} 
\]
and so 
\[
 \min\{b_1,\ldots, b_n\} \leq \frac{na^*}{|c|} 
\]
Because the sets $B_1,\ldots, B_n$ are discrete, the sets
\[
\tilde{B_i} =  \left\{b_i \in B_i:b_i \leq   \frac{na^*}{|c|}     \right\}
\]
are finite for $i=1,\ldots, n$, and so the set of fractions
\[
\mcf = \bigcup_{i=1}^n \left\{ \frac{\varepsilon_i a_i}{b_i} : a_i \in A_i, b_i \in \tilde{B_i}, \varepsilon_i \in \{1,-1\}    \right\}
\] 
is also finite.  
Every representation of $c$ in $\mce^{\pm}(\mca,\mcb)$ 
must include at least one fraction in the set \mcf.  
By the pigeonhole principle, 
if $r^{\pm}_{\mca,\mcb}(c) = \infty$, then there must exist $j_1 \in \{1,\ldots, n\}$ 
such that the fraction 
$\varepsilon_{j_1} a_{j_1}/b_{j_1} \in \mcf$ 
occurs in infinitely many representations.  
Let $j_1$ be the smallest integer in $\{1,\ldots, n\}$ with this property, 
and let $J_1 = (j_1)$.  
Let $L_1$ be the $(n-1)$-tuple obtained by deleting $j_1$ from $(1,\ldots, n)$, 
that is, 
\[
L_1 = (\ell_1,\ldots, \ell_{n-1}) = (1,\ldots, n) \setminus J_1.
\] 
We obtain 
\[
c_1 = c - \frac{\varepsilon_{j_1} a_{j_1}}{b_{j_1}} \in \mce^{\pm}(\mca_{L_1},\mcb_{L_1}) 
\]
and 
\[
r^{\pm}_{\mca_{L_1},\mcb_{L_1}}\left( c_1\right) 
= r^{\pm}_{\mca_{L_1},\mcb_{L_1}}\left( c - \frac{\varepsilon_{j_1} a_{j_1}}{b_{j_1}} \right) 
= \infty
\]
If $c_1 = 0$, then 
\[
c = \frac{\varepsilon_{j_1} a_{j_1}}{b_{j_1}} \in \mce^{\pm}(\mca_{J_1},\mcb_{J_1}) 
\subseteq \mcj(\mca,\mcb).
\]

If $c_1 \neq 0$, then we repeat this procedure.       
Because $r^{\pm}_{\mca_{L_1},\mcb_{L_1}}\left( c_1\right) = \infty$, 
we obtain $j_2 \in \{1,\ldots, n\}$ with $j_2 > j_1$ 
and a fraction $\varepsilon_{j_2} a_{j_2}/b_{j_2} \in \mcf$ 
that occurs in infinitely many representations of $c_1$.  
Let $j_2$ be the smallest integer in $\{j_1+1,\ldots, n\}$ with this property.  
Let $L_2$ be the $(n-2)$-tuple obtained by deleting $j_1$ and $j_2$ from $(1,\ldots, n)$, 
that is, $J_2 = (j_1, j_2)$ and $L_2 = (1,\ldots, n)\setminus  J_2$. 
Let 
\[
c_2 = c_1 - \frac{\varepsilon_{j_2} a_{j_2}}{b_{j_2} } 
= c  - \left(\frac{\varepsilon_{j_1} a_{j_1}}{b_{j_1} } 
+ \frac{\varepsilon_{j_2} a_{j_2}}{b_{j_2} } \right) 
\in \mce^{\pm}(\mca_{L_2},\mcb_{L_2}).
\]
We have proved that 
\[
r^{\pm}_{\mca_{L_2},\mcb_{L_2}}( c_2 ) = 
r^{\pm}_{\mca_{L_2},\mcb_{L_2}}\left( c - \left( \frac{\varepsilon_{j_1} a_{j_1}}{b_{j_1}} 
+ \frac{\varepsilon_{j_2} a_{j_2}}{b_{j_2}}\right)  \right) = \infty.
\]
If $c_2 = 0$, then 
\[
c = \frac{\varepsilon_{j_1} a_{j_1}}{b_{j_1}} 
+ \frac{\varepsilon_{j_2} a_{j_2}}{b_{j_2}} \in \mce^{\pm}(\mca_{J_2},\mcb_{J_2}) 
\subseteq \mcj(\mca,\mcb).
\]
If $c_2 \neq 0$, then we repeat this procedure.  

After $s$ iterations, we obtain the $s$-tuple 
\[
J_s = (j_1,\ldots, j_s) \preceq (1,\ldots, n),
\]
the $(n-s)$-tuple 
\[
L_s = (1,\ldots, n)\setminus (j_1,\ldots, j_s), 
\]
and fractions $\varepsilon_{j_i} a_{j_i}/b_{j_i}$ for $i =1,\ldots, s$ such that 
the weighted Egyptian number
\[
c_s = c - \sum_{i=1}^s \frac{ \varepsilon_{j_i} a_{j_i} }{ b_{j_i} } 
\in \mce^{\pm}(\mca_{L_s}, \mcb_{L_s})
\]
satisfies
\[
r^{\pm}_{\mca_{L_s},\mcb_{L_s}}  \left( c_s \right) = \infty.
\]
By Theorem~\ref{sierpinski:theorem:reps1-2}, this is impossible if $n-s = 2$ and $c_s \neq 0$.
Therefore, if $r^{\pm}_{\mca,\mcb}(c) = \infty$, then $c_s = 0$ for some 
$s \in \{1,\ldots, n-2\}$ and so 
\[
c \in  \mce^{\pm}(\mca_{J_s},\mcb_{J_s}) \subseteq \mcj(\mca,\mcb).
\]
This completes the proof. 
\end{proof}

\bl             \label{sierpinski:lemma:nowhereDense}
Let $A$ be a nonempty finite set of positive real numbers and let $B$ be an infinite 
discrete set of positive real numbers.  If $X$ is a nowhere dense set of real numbers,
then 
\[
Y = \left\{  x + \frac{\varepsilon a}{b}: x \in X, a \in A,  b \in B, \text{ and } 
 \varepsilon \in \{1,-1\}  \right\}
\]
is also a nowhere dense set of real numbers.     
\el

\begin{proof}
Let $a^* = \max(A)$.  
Because $X$ is nowhere dense, 
for every open interval $(u',v')$ there is a nonempty subinterval $(u,v)$ contained in in $(u',v')$ 
such that $X \cap (u,v) = \emptyset$.     Let $0 < \delta < (v-u)/2$ and let 
\[
y  \in \left(u+\delta, v-\delta \right) 
\]
for some 
\[
y = x + \frac{\varepsilon a}{b} \in Y.
\]
If $\varepsilon = 1$ and $y = x+a/b$, then 
\[
 x < y  < v - \delta < v
\] 
and so $x \leq u$.  Therefore, 
\[
x \leq u < u+\delta < y = x + \frac{a}{b} 
\]
and so $\delta < a/b$.
If $\varepsilon = -1$ and $y = x-a/b$, then 
\[
u <  u + \delta < y  < x
\] 
and so $x \geq v$.  Therefore, 
\[
x - \frac{a}{b} = y < v - \delta < v \leq x
\]
and $\delta < a/b$.
In both cases, $b < a/\delta \leq a^*/\delta$.
Because $A$ is finite and $B$ is discrete, the set 
\[
K = \left\{ 
\frac{\varepsilon a}{b}: a \in A,  b \in B , \varepsilon \in \{1,-1\}, 
\text{ and } b \leq \frac{a^*}{\delta}    \right\}
\]
is finite.   
We have  
\[
Z  =  \left\{ x + \kappa : x \in X \text{ and } \kappa \in K \  \right\} 
= X +K  \subseteq Y
\]
and 
\[
Y \cap (u+\delta, v- \delta) = Z \cap (u+\delta, v- \delta).
\]
The set $Z$ is the union of a finite number of translates of 
the nowhere dense set $X$.  
Because a translate of a nowhere dense set is nowhere dense, 
and because a finite union of nowhere dense sets is nowhere dense, 
it follows that $Z$ is nowhere dense.  
Therefore, the interval $(u+\delta, v- \delta)$ contains a nonempty open subinterval 
that is disjoint from $Y$. 
This completes the proof.
\end{proof}

\bt
Let $\mca = (A_1,\ldots, A_n)$ be a  sequence of nonempty finite sets 
of positive real numbers, 
and let $\mcb = (B_1,\ldots, B_n)$ be a sequence of infinite discrete 
sets of positive real numbers.  
The set $\mce^{\pm}(\mca,\mcb)$ is nowhere dense.
\et

\begin{proof}
The proof is by induction on $n$.  
If $n=1$, then  $\mca = (A_1)$,  $\mcb = (B_1)$, 
the set $\mce^{\pm}(\mca,\mcb)$ is discrete, and a discrete set is nowhere dense.
The inductive step follows immediately from Lemma~\ref{sierpinski:lemma:nowhereDense}.
\end{proof}

\emph{Acknowledgements.}  David Ross introduced me to this subject 
at the AIM workshop ``Nonstandard methods in  combinatorial number theory.''

This paper was written on Nantucket in August, 2017.  
I thank The Bean and the Nantucket Atheneum 
for providing excellent work environments.

\def\cprime{$'$} \def\cprime{$'$}
\providecommand{\bysame}{\leavevmode\hbox to3em{\hrulefill}\thinspace}
\providecommand{\MR}{\relax\ifhmode\unskip\space\fi MR }
\providecommand{\MRhref}[2]{%
  \href{http://www.ams.org/mathscinet-getitem?mr=#1}{#2}
}
\providecommand{\href}[2]{#2}

\end{document}